\documentclass[reqno,12pt]{article}

\usepackage{a4wide}
\usepackage{amsmath}
\usepackage{amssymb}
\usepackage{amsthm}
\usepackage[utf8]{inputenc} 
\usepackage{graphicx} 
\usepackage[english, russian]{babel}

\numberwithin{equation}{section}

\newtheorem{theo}{Theorem}
\newtheorem{conj}{Conjecture}
\newtheorem{coro}{Corollary}

\theoremstyle{remark}

\newcommand{\Z}{\mathbb{Z}}
\newcommand{\Q}{\mathbb{Q}}

\renewcommand{\C}{\mathbb{C}}
\newcommand{\K}{\mathbb{K}}

\newcommand{\Qbar}{\overline{\mathbb Q}}
\newcommand{\Qb}{\overline{\mathbb Q}}

\newcommand{\etoile}{^\star}

\newcommand{\unmu}{\{1,\ldots,\mu\}}

\newcommand{\Li}{{\rm Li}}

\newcommand{\bz}{b_0}
\newcommand{\dd}{{\rm d}}
\newcommand{\tra}{  ^t \! \!}
\newcommand{\GLmu}{{\rm GL}_\mu}
\newcommand{\AH}{A_{[H]}}
\newcommand{\calY}{{\mathcal Y}}
\newcommand{\calN}{{\mathcal N}}
\newcommand{\OK}{{\mathcal O}_{\K}}
\newcommand{\PK}{{\mathcal P}_{\K}}
\newcommand{\pe}[1]{\left\lfloor #1 \right\rfloor}
\newcommand{\pes}[1]{\left\lceil #1 \right\rceil}

\newcommand{\calG}{{\mathcal G}}
\newcommand{\calF}{{\mathcal F}}

\begin{document}

 \selectlanguage{english}

\title{On the denominators of the Taylor coefficients of $G$-functions}
\date\today
\author{S. Fischler and T. Rivoal}
\maketitle

\begin{abstract}  Let $\sum_{n=0}^\infty a_n z^n\in \Qbar[[z]]$ be a $G$-function, 
and, for any $n\ge0$, let $\delta_n\ge 1$ denote the least integer such that $\delta_n a_0, \delta_n a_1, \ldots, \delta_n a_n$ are all algebraic integers. By definition of a $G$-function, there exists some constant $c\ge 1$ such that $\delta_n\le c^{n+1}$ for all $n\ge 0$. In practice, it is observed that $\delta_n$  always divides
$D_{bn}^{s} C^{n+1}$ where $D_n=\textup{lcm}\{1,2, \ldots, n\}$,  $b, C$ are positive integers and $s\ge 0$ is an integer.  
We prove that this observation holds for any $G$-function  provided the following conjecture is assumed: {\em Let  $\K$ be a number field, and $L\in\K[z,\frac{\dd }{\dd z}]$ be a $G$-operator; then the generic radius of solvability $R_v(L)$  is equal to 1, for all finite places $v$ of $\K$ except a finite number.} The proof makes use of very precise estimates in the theory of $p$-adic differential equations, in particular the Christol-Dwork Theorem. Our result becomes unconditional when $L$ is a geometric differential operator, a special type of $G$-operators for which the  conjecture is known to be true. The famous Bombieri-Dwork Conjecture asserts that any $G$-operator is of geometric type, hence it implies the above conjecture.
\end{abstract}

\section{Introduction}

Siegel~\cite{siegel} introduced  
the class of $G$-functions as a vast generalization of the geometric series $\sum_{n=0}^\infty z^n=\frac{1}{1-z}$ and its primitive  
$-\log(1-z)=\sum_{n=1}^\infty \frac{z^n}{n}$. 
A $G$-function is a formal power series $g(z)=\sum_{n=0}^{\infty} a_n z^n\in \Qbar[[z]]$ 
such that there exists $C>0$ such that (upon fixing an embedding of 
$\Qbar$ into $\C$):

$(i)$ the maximum of the moduli of the conjugates of $a_n$ is $\leq C^{n+1}$.

$(ii)$ there exists a sequence of integers $\delta_n$ such that $0 <\delta_n  \leq C^{n+1}$ and 
$\delta_na_m$ is an algebraic integer for all~$m\le n$.

$(iii)$ $f(z)$ is holonomic, i.e. it satisfies an homogeneous linear differential equation with 
coefficients in $\Qb(z)$.

\noindent These properties ensure  that $g(z)$ is holomorphic at $z=0$.
Prototypical examples of $G$-functions are the algebraic functions over $\Qbar(z)$   holomorphic at $z=0$, or the hypergeometric series ${}_{p+1}F_p(z)$ with rational parameters. The latter class includes the polylogarithmic functions  $\Li_k(z)=\sum_{n=1}^\infty \frac{z^n}{n^k}$, where $k$ is any given integer. Any holonomic power series with integer coefficients and positive radius of convergence is a $G$-function.

\medskip

In this paper, we are interested in the sequence of denominators $\delta_n$ in $(ii)$. In practice, it has been observed for long that the least denominator $\delta_n$  always divides
$D_{bn}^{s} C^{n+1}$ where $D_n=\textup{lcm}\{1,2, \ldots, n\}$,  $b, C$ are positive integers and $s\ge 0$ is an integer. We recall that $D_n= e^{n(1+o(1))}$ by the prime number theorem, so that $D_{bn}^{s} C^{n+1}$ has geometric growth. By Eisenstein's theorem, we have $s=0$ for any algebraic function  over $\Qbar(z)$. But this is not the general situation, for we have $s=k$ in the case of the polylogarithms $\Li_k(z)$ for $k\ge 1$, because $\textup{lcm}\{1^k,2^k, \ldots, n^k\}=\textup{lcm}\{1,2, \ldots, n\}^k$.

Our main result is a proof of this observation, in very precise terms, conditionally on a well-known conjecture on the generic $v$-adic radii of convergence of solutions of $G$-operators, but also unconditionally for an important class of $G$-operators. We recall that a differential operator $L\in \Qbar[z, \frac{d}{dz}]$ is said to be a $G$-operator if the differential equation $Ly(z)=0$ has a $G$-function  solution for which it is minimal for the degree (or order) in $\frac{d}{dz}$. There is no loss of generality in assuming from the beginning that $L\in \K[z, \frac{d}{dz}]$ for some  number field $\K$. Moreover, we can trivially rewrite the differential equation $Ly(z)=0$, of order $\mu$ say, as a differential system  $X'(z)=G(z)X(z)$ of order $1$, where $X(z)=\tra \big(y(z), \ldots, y^{(\mu-1)}(z)\big)$ and $G(z)$ is a square matrix in $M_{\mu}(\K(z))$. For any finite place $v$ of $\K$, the generic radius of solvability $R_v(L)$ of $L$ at $v$ is the greatest real number $R \leq 1$ such that the system $X'(z)=G(z)X(z)$ has a solution matrix ${\mathcal X}(z)$ at the generic point $t$  relative to $v$ consisting in functions analytic in $D(t, R^-)$. See \cite[pp. 67 and 70]{AndreGgeom} or \cite{DGS} for these classical definitions.  
The following conjecture is often attributed to Bombieri.

\begin{conj} \label{conj1} Let  $\K$ be a number field, and $L\in\K[z,\frac{\dd }{\dd z}]$ be a $G$-operator. Then the generic radius of solvability $R_v(L)$  is equal to 1, for all finite places $v$ of $\K$ except a finite number.
\end{conj}

This conjecture is a refinement of the following result, 
 which is a consequence of two independent theorems, of Andr\'e and Bombieri~\cite[p. 228, Theorem 2.1]{DGS} on the one hand, and of the Chudnovsky's~\cite[p. 17, Theorem III]{chud} on the other hand: if $L\in\K[z,\frac{\dd }{\dd z}]$ is a $G$-operator, then $\sum_{v} \log(1/R_v(L))<\infty$ where the summation is over the finite places $v$ of $\K$; see~\cite[p. 228, Theorem 2.1]{DGS}.
The conjecture is in fact  known to be true when $L$ is a $G$-operator of a special type, namely when it is a {\em geometric differential operator}, i.e. a product of factors of Picard-Fuchs differential operators; see \cite[pp. 39 and 111]{AndreGgeom}. Roughly speaking, a Picard-Fuchs equation is the differential equation satified by the periods of a family of algebraic varieties (defined over $\Qbar$ in our situation) depending on a complex parameter.  
It follows that Conjecture \ref{conj1} is a consequence of the Bombieri-Dwork conjecture, as formulated by Andr\'e~\cite[p. 111]{AndreGgeom}, which asserts the converse: every $G$-operator is a geometric differential operator, a conjecture usually phrased  in looser terms as ``every $G$-operator comes from geometry''. 

\medskip

Our main result is the following. In the case where 0 is not a singularity of $L$, it is implicit in Bombieri's paper \cite[\S 6]{Bombieri}; see also \cite[p.~74, Theorem 5.2]{AndreGgeom} and \S 2 below.

\begin{theo} \label{thprinc}
Let $g(z) = \sum_{n=0}^\infty a_n z^n $ be a $G$-function annihilated by a $G$-operator $L$ of order $\mu$ for which Conjecture \ref{conj1} holds. Then there exist some positive integers $b$, $\bz$, $C$ such that, for any $n\geq 0$, $D_{bn+\bz}^{\mu-1} C^{n+1} a_n$ is an algebraic integer.

Moreover $b$ can be chosen as the least common multiple of the denominators of the exponents of $L$ at 0.
\end{theo}

We do not need to assume that $L$ is the {\em minimal} $G$-operator that annihilates $g(z)$, though this assumption obviously provides the smallest value of $\mu$ we can achieve by our method. The exponent $\mu-1$ can not be improved in general as the example of polylogarithms shows, but algebraic functions also show it is not always best possible. In the final section, we provide  another example where the exponent $\mu-1$ is optimal. If $0$ is not a singularity of $L$, then we can take $b=1.$

From the above comments, we also deduce:

\begin{theo} \label{coro1}
The conclusions of Theorem \ref{thprinc} hold unconditionally if $L$  is a geometric differential operator.
\end{theo}

Since any $G$-function is annihilated by a $G$-operator, we also deduce the following corollary to Theorem~\ref{thprinc} (because $D_{bn+\bz}^s C^{n+1} $ divides $D_{(b+\bz)n}^s (D_{\bz}^s C)^{n+1}$ for any $n\ge 0$). 
Again, it is unconditional if the $G$-operator is a geometric one.

\begin{coro}  \label{coro2}
Assume that Conjecture \ref{conj1} holds. Then for any $G$-function $\sum_{n=0}^\infty a_n z^n $  there exist some positive integers $b$,    $C$ and an integer $s\ge 0$ such that, for any $n\geq 0$, $D_{bn }^{s} C^{n+1} a_n$ is an algebraic integer.
\end{coro}

As an application of this result, let us define a filtration on the set $\calG$ of all $G$-functions as follows. Given $s\geq 0$, we denote by $\calG_s$ the set of $g(z) = \sum_{n=0}^\infty a_n z^n  \in\calG$ with the following property: there exist $b,  C\geq 1$ such that, for any $n\geq 0$, $D_{bn }^{s} C^{n+1} a_n$ is an algebraic integer. Then $\calG_s$ is a $\Qbar$-vector space stable under derivation, and we have
$$\calG_s\subset\calG_{s+1} \;\mbox{ and }\; \calG_s\calG_t\subset\calG_{s+t} \;\mbox{ for any }\; s,t\geq 0.$$
Elements of $\calG_0$ are exactly globally bounded $G$-functions (see  \cite{christol}); for any $s\geq 1$ the $s$-th polylogarithm $\Li_s(z)  $ belongs to $\calG_s$ and not to $\calG_{s-1}$. Given a hypergeometric function $g(z) = {}_{p+1}F_p(z)$ with rational coefficients, Christol's criterion~\cite{christol} (see also \cite{drr}) enables one to compute the least integer $s$ such that $g\in\calG_s$. Letting $\calG_{\infty} = \cup_{s=0}^\infty \calG_s$,  Corollary \ref{coro2} shows that  Conjecture \ref{conj1}  implies $\calG_\infty = \calG$ (i.e., $(\calG_s)_{s\geq 0}$ is a total filtration of $\calG$). As the proof of Theorem \ref{thprinc} shows (see \S \ref{subsec34}), $\calG_s$ is the set of $G$-functions $ g(z) = \sum_{n=0}^\infty a_n z^n$ with the following property (where $\K$ is a number field containing all coefficients $a_n$): there exists $b \geq 0$ such that for any sufficiently large prime number $p$, for any place $v$ of $\K$ above $p$, and for any $n\geq 0$, we have $v(D_{bn  }^s a_n)\geq 0$. 

It would be very interesting to find a graduation $\calG_\infty = \oplus_{t=0}^\infty \calF_t$ such that $\calG_s = \oplus_{t=0}^s  \calF_t$ for any $s\geq 0$ and $\calF_t\calF_{t'}\subset\calF_{t+t'}$ for any $t,t'\geq 0$.

\bigskip

The structure of this text is as follows. In \S \ref{sec2} we comment on some applications of our results, and relate them to linear recurrences with polynomial coefficients. Then we prove Theorem \ref{thprinc} in \S \ref{sec3}, and apply it in \S \ref{sec4} in a situation related to Ap\'ery's proof that $\zeta(3)$ is irrational. 

\bigskip

\noindent {\bf Acknowledgment.}  We warmly thank Yves Andr\'e and Gilles Christol for many discussions on the subject of this paper and for patiently answering our questions. 

\section{Some comments} \label{sec2}

When a $G$-function is a solution of a $G$-operator 
{\em regular} at $z=0$, 
our proof of  Theorem~\ref{thprinc}  can be simplified 
and appears already implicitly in \cite[\S 6]{Bombieri} (see also \cite[p.~74, Theorem 5.2]{AndreGgeom}). In this case Theorem~\ref{thprinc} follows from    the very general theorem of Dwork-Robba~\cite[p. 120]{DGS}; in contrary to the proof below, all finite places can be dealt with in the same way, and no shearing transformation is necessary.

However, not every $G$-function can be solution of a $G$-operator regular at $z=0$. Indeed, assume that an operator $L\in \Qbar[z,\frac{d}{dz}]$
is minimal for some $G$-function $g(z)$  ($L$ is then a $G$-operator) and that another solution $h(z)$ at $z=0$ of $Ly(z)=0$ has a non-trivial monodromy at $z=0$. (For instance, we can consider the hypergeometric function $g(z)={}_2F_1[\frac12,\frac12;1;z]$ with a second solution $h(z)=f(z)+\log(z)g(z)$ for some $f(z)\in\mathbb Q[[z]]$.)  Since $L$ is a right factor of any differential operator $M\in \Qbar[z,\frac d{dz}]$ such that $Mg(z)=0$, we have $Mh(z)=0$ as well, so that $z=0$ is a singularity of $M.$

Therefore we need a generalization of the Dwork-Robba Theorem. A first step in this direction has been made by 
 Adolphson,  Dwork and   Sperber \cite{ADS}: their result covers the case of a single Jordan block, and can be used to derive Theorem~\ref{thprinc} in this case. The general theorem we shall use is due to 
 Christol-Dwork~\cite{Duke}; we will make a certain number of transformations to be in position to apply it.

\bigskip

So far, there is no known intrinsic way to check if a differential operator $L$ is a $G$-operator, except if we can find a $G$-function it annihilates and for which it is minimal, or if we know that it comes from geometry.~(\footnote{Another conjecture predicts that  a differential operator in $\mathbb Q[z,\frac d{dz}]$ is a $G$-operator if and only if its $p$-curvature is nilpotent for all but finitely many primes $p$. (There is a similar conjecture when $\mathbb Q$ is replaced by any number field.) The ``only if'' part is a consequence of a theorem of the Chudnovsky's~\cite{chud}. 
The $p$-curvature can be computed (see~\cite{bostan2} for an efficient algorithm) from the associated differential system for any $p$ without computing an exact solution of the differential equation. But the drawback of the conjecture is that it has to be computed for {\em infinitely} many $p$, which makes it unpractical. However, if the $p$-curvature is observed to be nilpotent  for a very large number of primes $p$, the conjecture can be used as an oracle to guess if an operator is a $G$-operator (and then try to prove this by finding an explicit $G$-function solution). This procedure was used with success in~\cite{bostan}.}) Moreover, there is no known algorithm to find such a $G$-function solution, and in 
practice the situation is reversed: we are given an explicit $G$-function $g(z)$ and  then  we compute its associated $G$-operator $L$. Thus, one may think that Theorems~\ref{thprinc} and~\ref{coro1} are of limited practical interest because we can directly see that the denominators of the Taylor coefficients of $g(z)$ are of the desired form. 

However, they do have a practical interest, especially if $L$ is a geometric operator because the results are then unconditional. Indeed, let us assume we are given  an explicit $G$-function, or a period function, and its associated $G$-operator $L$. By a well-known theorem of Andr\'e, Chudnovsky and Katz (see~\cite{andre} or~\cite{firi} for a precise statement), we know that at any $\alpha\in \Qbar \cup\{\infty\}$, a basis of local solutions at $z=\alpha$ of $Ly(z)=0$ can be made of 
\begin{itemize}
\item[(1)] $G$-functions in the variable $z-\alpha$, if $\alpha$ is a regular point of $L$;
\item[(2)] linear combinations of $G$-functions  in the variable $z-\alpha$, with coefficients some polynomials (over $\Qbar$) in certain rational powers of $z-\alpha$ and integral powers of $\log(z-\alpha)$, if $\alpha$ is a singularity of $L$. 
\end{itemize}
(Above and below, we replace $z-\alpha$ by $1/z$ if $\alpha=\infty$.) In general, we know nothing about the $G$-functions appearing in these basis, in the sense that there is no general procedure to find  closed form formulas for their Taylor coefficients. But Corollary \ref{coro2} applies to these 
$G$-functions,  which thus provides non-trivial informations on the denominators of their coefficients. In case (1), the $G$-functions are solutions of $Ly(z)=0$, hence we can take $s=\mu-1$. In case (2), the $G$-functions are no longer necessarily solutions of $Ly(z)=0$ but they are solutions of $G$-operators of order at most $\mu^2$ (see \cite[p. 110]{andre}) so that we can take $s=\mu^2-1$.
This corollary may even be applied  unconditionally because if $L$ is a geometric operator, this is also the case of any of the shifted operators $L_\alpha$ by any algebraic number $\alpha$, where by definition $L_\alpha y(z)=0$ is equivalent to $Ly(z-\alpha)=0$. 

\medskip

The above remarks can be translated in the language of linear recurrences with polynomial coefficients. We consider a sequence $(a_n)_{n\ge 0}$ of algebraic numbers satisfying a recurrence of the form $\sum_{j=0}^k p_j(n) a_{n+j}=0$, where $p_j(n)\in \Qbar[n]$. 
Such recurrences naturally occur (sometimes in many different ways) when trying to prove  the irrationality of constants like $\pi, \log(2), \zeta(3)$, and it is important to control the growth of the sequence of denominators $\delta_n$ of $a_n$. 
The general situation is that $\delta_n$ grows like $(b n+\bz)!^s$, and not that $\delta_n$ divides $D_{bn+\bz}^s C^{n+1}$. However, if it is of the latter type, then the series $g(z)=\sum_{n=0}^\infty a_n z^n$ is not far from being a $G$-function because the recurrence relation  gives rise to a linear differential equation $Ly(z)=0$ with coefficients in $\Qbar(z)$ of which $g(z)$ is solution, and indeed $g(z)$ is a $G$-function if $a_n\in \mathbb Q$ for instance. 

Consider a $G$-operator $L$ and assume that there exists a local solution $\sum_{n=0}^\infty v_n (z-\alpha)^n$ of $Ly(z)=0$ at some given $\alpha \in \Qbar \cup\{\infty\}$. To $L$, we can associate bijectively a linear recurrence  
\begin{equation}\label{eq:rec}
\sum_{j=0}^k q_{j}(n) V_{n+j}=0
\end{equation}
 with $q_{j}(n)\in \Qbar[n]$ of which $(v_n)_{n\ge 0}$ 
is a solution.

\begin{coro}\label{coro3} If Conjecture~\ref{conj1} holds, then for any solution $(V_n)_{n\ge 0}$ of the recurrence~\eqref{eq:rec} with $V_n\in \Qbar$  for all $n\ge 0$, we can find some positive integers $b,  C$ and an integer $s\ge 0$ such that $D_{bn }^s C^{n+1} V_n$ is an algebraic integer  for all $n\ge 0$. This is unconditional if $L$ is a geometric operator.  
\end{coro}
Indeed, any solution $(V_n)_{n\ge0}$ of~\eqref{eq:rec} is such that $\sum_{n=0}^\infty V_n (z-\alpha)^n$ is a solution of $Ly(z)=P(z)$ for some polynomial $P(z)\in \Qbar[z]$ of degree $k$ say ($P$ arises because of the initial conditions). Then 
$(\frac{d}{dz})^{k+1} \circ L $ is a $G$-operator that annihilates $\sum_{n=0}^\infty V_n (z-\alpha)^n$. 
Hence, $\sum_{n=0}^\infty V_n z^n$ is a $G$-function by the Andr\'e-Chudnovsky-Katz theorem. We conclude by an application of Corollary~\ref{coro2}. The problem of bounding the growth of denominators of rational solutions of Ap\'ery like linear recurrences in absence of explicit formulas is an important one in Diophantine approximation. It was considered in~\cite{beukers11} and~\cite{zagier} for instance.

In practice, we start with an explicit $G$-function $\sum_{n=0}^\infty a_n z^n$ and compute its minimal $G$-operator. The above property  enables us to construct new linear recurrences with solutions with controlled denominators without having to construct explicitely these solutions. In fact, we don't how to do that in general, except when $\alpha=0$ and $V_n=a_n$ is the sequence we start with. If the $G$-operator comes from geometry, this is unconditional. We develop this procedure in  Section~\ref{sec:example} on an explicit example.

\section{Proof of Theorem~\ref{thprinc}} \label{sec3}

We split the proof of Theorem \ref{thprinc} in several parts. The strategy is to apply the Christol-Dwork Theorem at all (but finitely many) finite places. To meet the nilpotence assumption $(ii)$ in Theorem \ref{theoCD} below, we use two preliminary steps (\S\S \ref{subsec31} and \ref{subsec32}). Then we recall and apply the Christol-Dwork Theorem in \S \ref{subsec33}, and conclude the proof in \S \ref{subsec34} by considering the (finitely many) remaining finite places.
The constant $b$ (resp. $\bz$ and $C$) is   constructed in \S \ref{subsec31} (resp. after Eq. \eqref{eqgun} and in \S \ref{subsec34}).

Unless otherwise stated, all the tools we  use appear in \cite{DGS}.

\subsection{Reduction to the case of integer exponents at 0} \label{subsec31}

We first recall that from the Andr\'e-Chudnovsky-Katz Theorem, we know that a $G$-operator is Fuchsian with rational exponents at every algebraic point (and at  infinity).

We shall prove  Theorem \ref{thprinc} in a special case below (\S\S \ref{subsec32} -- \ref{subsec34}), namely when all exponents of $L$ at 0 are integers (so that $b=1$). To deduce the general case, we recall that all exponents of $L$ at 0 are rational numbers, and denote by $b\geq 1$  the least common multiple of their denominators. Then we perform the change of variable $x = z^{1/b}$, by letting $\widetilde g (x) = \sum_{n=0}^\infty a_n x^{bn}$. Let $\widetilde L \in \Qbar[x, \frac{\dd}{\dd x}]$ be obtained from $L$ in this way, so that $\widetilde L \widetilde g(x) = 0$. The exponents of $\widetilde L $ at 0 are those of $L$ multiplied by $b$: they are integers. Moreover Conjecture~\ref{conj1} holds for $\widetilde L$, and  $\widetilde L$ has order $\mu$: the special case of  Theorem \ref{thprinc} we shall prove below yields $\bz$ and  $C$ such that   $D_{bn+\bz}^{\mu-1} C^{bn+b} a_n$ is an algebraic integer  for any $n\geq 0$. This concludes the proof of the general case.

\subsection{Shearing transformation} \label{subsec32}

From now on, we shall prove  Theorem \ref{thprinc} under the additional assumption that all exponents of $L$ at 0 are integers (so that $b=1$).
We write $L = B_\mu(z) (  \frac{\dd}{\dd z})^\mu + \ldots + B_0(z)$ with $B_0(z),\ldots, B_\mu(z) \in\Qbar[z]$ and $B_\mu(z)\neq 0$. 
To $L$, we associate the differential system $zY'(z)=A(z)Y(z)$, with
$$A(z) = 
\left( \begin{matrix}
0&1&0&0&\ldots &0&0\\
0&1&1&0&\ldots &0&0\\
0&0&2&1&\ldots &0&0\\
0&0&0&3&\ldots &0&0\\
\vdots & \vdots & \vdots & \vdots & & \vdots & \vdots \\
0&0&0&0&\ldots &\mu-2 &1 \\
-a_0&-a_1&-a_2&-a_3&\ldots &-a_{\mu-2}&\mu-1-a_{\mu-1}
\end{matrix}\right)
\in M_\mu(\Qbar(z))$$
where $a_i(z) = z^{\mu-i} B_i(z) / B_\mu(z)$ is holomorphic at $z=0$ for $i \in \{0,\ldots,\mu-1\}$ since $0$ is a regular singularity of $L$. 
We have $Ly(z)=0$ if, and only if, $zY'(z) =A(z)Y(z)$ where $Y(z) = \tra \big(y(z),zy(z)',z^2y''(z),\ldots,z^{\mu-1}y^{(\mu-1)}(z)\big)$; 
see for instance \cite[Chapter 6]{Coddington-Carlson}.

This is not exactly the usual way to write the differential system $X'(z)=G(z)X(z)$ associated to $L$ for the vector $X(z) = \tra \big(y(z),y(z)',y''(z),\ldots,y^{(\mu-1)}(z)\big)$, but this is a standard variation when dealing with regular singularities. In particular, the generic radii of solvability of $L$ can be computed  with the system $X'(z)=G(z)X(z)$ or with the system $zX'(z)=A(z)X(z)$, giving the same results.

\medskip

The eigenvalues of $A(0)$ are exactly the exponents of $L$ at 0: by assumption they are integers. Therefore they can be made all equal to 0 by applying a suitable shearing transformation. In precise terms, for any $H(z) \in \GLmu (\Qbar[z,z^{-1}])$ (i.e., such that both $H(z)$ and $H(z)^{-1}$ have entries in $\Qbar[z,z^{-1}]$), let 
$$\AH(z) =  zH'(z)  H(z)^{-1} + H(z) A(z) H(z)^{-1}.$$
Then upon letting $Z(z) = H(z)Y(z)$ we have $zY'(z)=A(z)Y(z)$ if, and only if, $zZ'(z) = \AH(z) Z(z)$. Using the matrix $H(z)$ provided by \cite[p. 106, Corollary 8.3]{DGS} we have the additional properties that $\AH(z)\in M_\mu(\Qbar[[z]])$ and $0$ is the only eigenvalue of $\AH(0)$, so that $\AH(0)$ is nilpotent.

\medskip

Since the differential system $zZ'(z) = \AH(z) Z(z)$ has a regular singularity at 0, it has a solution matrix of the form $U(z) z^{\AH(0)}$ with $U(z) \in \GLmu(\Qbar[[z]])$ such that $U(0) = I_\mu$. Then $\calY (z) = H(z)^{-1} U(z) z^{\AH(0)}$ is a  solution matrix of the differential system $zY'(z)=A(z)Y(z)$. Since $Y_0 (z) = \tra (g(z),zg'(z),  \ldots,z^{\mu-1}g^{(\mu-1)}(z)) \in M_{q,1}(\Qbar[[z]])$ is a solution of this system, there exists $L_0 = \tra (\ell_1,\ldots,\ell_\mu)\in M_{q,1}(\Qbar)$ such that $Y_0(z) = \calY(z) L_0$, so that
\begin{equation} \label{eqgun}
g(z) = \ell_1 y_{1,1}(z) + \ldots +  \ell_\mu y_{1,\mu}(z) 
\end{equation}
where $y_{i,j}(z)$ is the coefficient of index $(i,j)$ of $\calY(z)$. We denote by $\sum_{t= -\bz}^{\bz} h_{i,j,t} z^t$ the coefficient of index $(i,j)$ of $H(z)^{-1}$ and by $\sum_{n=0}^\infty u_{i,j,n}z^n$ that of $U(z)$, where $\bz\geq 0$ is fixed and $h_{i,j,t}, u_{i,j,n}\in\Qbar$. Since $\AH(0)$ is nilpotent, the coefficient of index $(i,j)$ of $z^{\AH(0)}$ is a polynomial $P_{i,j}(X)\in\Qbar[X]$ taken at $X=\log(z)$. Therefore the equality $\calY (z) = H(z)^{-1} U(z) z^{\AH(0)}$ yields
$$y_{1,i}(z) = \sum_{j,k=1}^\mu \Big( \sum_{t= -\bz}^{\bz} h_{1,j,t} z^t \Big)  \Big(  \sum_{n=0}^\infty u_{j,k,n}z^n \Big)  P_{k,i}\big(\log(z)\big)$$
so that using Eq. \eqref{eqgun} we have 
$$g(z) =   \sum_{i,j,k=1}^\mu \ell_i P_{k,i}\big(\log(z)\big)  \sum_{m= -\bz}^\infty  \Big( \sum_{t=-\bz} ^{\min(\bz,m)} h_{1,j,t}  u_{j,k,m-t}  \Big)z^m.  $$
Recall that $g(z) = \sum_{n=0}^\infty a_n z^n \in\Qbar[[z]]$; then we have the following expression of $a_n$, in which the polynomials $P_{k,i}$ appear only through their constant terms:
 \begin{equation} \label{eqansomme}
 a_n =   \sum_{i,j,k=1}^\mu \sum_{t=-\bz} ^{\min(\bz,n)}  \ell_i P_{k,i}(0)    h_{1,j,t}  u_{j,k,n-t} 
\end{equation}
This expression allows us to focus now on the coefficients $u_{j,k,n-t}$.

\subsection{Application of the Christol-Dwork Theorem} \label{subsec33}

In this section, we apply  the Christol-Dwork Theorem~\cite{Duke} in the setting of \S \ref{subsec32}; we shall deduce Theorem \ref{thprinc} in \S \ref{subsec34} below.   We use the following particular case of the version stated in~\cite[p. 159, Theorem 2.1]{DGS}. 
\begin{theo}[Christol-Dwork] \label{theoCD} Let $\K$ be a number field and $v$ a finite   place  of $\K$ over a prime $p$.
Let $G(z)\in M_\mu(\K(z) )$ and consider the differential system $zZ'(z)=G(z)Z(z)$. Assume  that:
\begin{enumerate}
\item[$(i)$]  The coefficients of $G(z)$ have no pole in the disk $D(0, 1^-)$.
\item[$(ii)$] The matrix  $G(0)$ is nilpotent.
\item[$(iii)$] We have  $\vert G(z)\vert_{{\rm Gauss},v}\le 1$.
\item[$(iv)$]   There exists a matrix solution $U_{t}(z)$ at the generic point $t$ (normalized by $U_{t}(t)=I_\mu$) which converges in $D(t,1^-)$.
\end{enumerate}
Let $U(z)z^{G(0)}$ be a matrix solution of  $z Z'(z)=G(z)Z(z)$ with  $U(z)=\sum_{n=0}^\infty U_n z^n$ where $U_0=I_\mu$ and $U_n\in M_\mu(\K)$. Then for all $n\ge0$, we have 
\begin{equation}\label{eq:CD}
\vert U_n\vert_v \le T(n)^{\mu-1+v_p((\mu-1)!)+\beta_\mu},
\end{equation}
where $T(n)  = \max_{1\leq k \leq n} |k|_v^{-1}$ if $n\geq 1$, $T(0)=1$, $\beta_\mu = \min\big(\mu-1,v_p(\prod_{j=1}^\mu \binom{\mu}{j})\big)$,   $v_p$ is the usual $p$-adic valuation on $\Q\etoile$, and $|M|_v = \max _{i,j} |m_{i,j}|_v$ for $M = (m_{i,j})\in M_\mu(\K)$.
\end{theo}

To apply this result, we let
  $\K$ be a number field such that $H(z)$, $H(z)^{-1}$, $A(z)$, $\AH(z)$ belong to $M_\mu(\K(z))$ and the numbers $\ell_i$ of Eq. \eqref{eqgun} belong to $\K$. Let $\PK$ denote the set of all finite places $v$ of $\K$, i.e. ultrametric absolute values of $\K$ normalized by $|p|_v = 1/p$ if $v$ is above the prime number $p$.  We shall check the assumptions of Theorem \ref{theoCD} with $G(z)=\AH(z)$, provided $p$ is large enough. For any $v\in\PK$ above $p$, the completion $\K_v$ of $\K$ with respect to $v$ is a finite extension of $\Q_p$ and for any $\alpha\in\K$ we have $|\alpha|_v = |N_{\K_v/\Q_p}(\alpha)|_p^{1/d_v}$ where $d_v = [\K_v:\Q_p]$; here $|x|_p$ is the usual $p$-adic absolute value on $\Q$, with $|p|_{p} = 1/p$. 

\medskip

Let $v\in \PK$, and $t$ be a generic point with respect to $v$. Then the local solutions at $t$ of $Ly(z)=0$, $zY'(z) = A(z)Y(z)$ and $zZ'(z) = \AH(z)Z(z)$ are related to one another as above. Since by assumption the differential operator $L$ has generic radius of solvability   equal to 1 (except for finitely many places $v$ over some number field containing its coefficients), the same property holds over $\K$ and  for the differential system  $zZ'(z) = \AH(z)Z(z)$. Let $p_1$ be such that this  radius is equal to 1 (so that assumption $(iv)$ of Theorem \ref{theoCD} holds) for any $v$ above any prime $p \geq p_1$.

\medskip

For $i,j\in\unmu$, we denote by $\frac{R_{i,j}(z)}{S_{i,j}(z)}$ the coefficient of index $(i,j)$ of the matrix $\AH(z)$, with $R_{i,j}(z) , S_{i,j}(z) \in  \OK[z]$ and $S_{i,j}(0)\neq 0$; here $\OK$ is the ring of integers of $\K$. Let $\omega_{i,j,\ell}$, $1\leq \ell \leq \deg S_{i,j}$, denote the roots of  $S_{i,j}(z)$. Enlarging $\K$ if necessary, we may assume that  $\omega_{i,j,\ell}\in\K$ for any $i$, $j$, $\ell$; moreover $\omega_{i,j,\ell}\neq 0 $ since $S_{i,j}(0)\neq 0$. At last, let $r_{i,j,t}$ and $s_{i,j,t}$ denote the non-zero coefficients of $R_{i,j}(z)$ and $S_{i,j}(z)$.  Let $\calN_2\in\K\etoile$ denote the product of all numbers $\omega_{i,j,\ell}$, $r_{i,j,t}$ and $s_{i,j,t}$ as $i$, $j$, $\ell$, $t$ vary. Then there exists $p_2$ such that for any prime number  $p \geq p_2$, the $p$-adic valuation of $N_{\K /\Q }(\calN_2) \in \Q\etoile$ is 0.

\medskip

At last, let $\calN_3\in\K\etoile$  denote the product of all numbers $\ell_i$, $P_{k,i}(0)$ and $h_{1,j,t}$ of Eq. \eqref{eqansomme} which are non-zero, as $i$, $j$, $k$, $t$ vary; then  there exists $p_3$ such that for any prime   $p \geq p_3$, the $p$-adic valuation of $N_{\K /\Q }(\calN_3) \in \Q\etoile$ is 0.

\medskip

Let $p_0 = \max(\mu + 1,p_1,p_2,p_3)$,  $p$ be a prime number such that $p\geq p_0$, and $v\in\PK$ be a place above $p$. 
Since $p\geq p_2$ we have $| \omega_{i,j,\ell}| _v = | r_{i,j,t}|_v = |s_{i,j,t}|_v = 1$ for any  $i$, $j$, $\ell$, $t$, so that $|R_{i,j}(z)|_{{\rm Gauss}, v} = \max_t  | r_{i,j,t}|_v = 1$,  $|S_{i,j}(z)|_{{\rm Gauss}, v} = 1$, and $|\AH(z)|_{{\rm Gauss}, v} = 1$. 
Moreover, the coefficients $\frac{R_{i,j}(z)}{S_{i,j}(z)}$ of $\AH(z)$ have no pole in the disk $D(0, 1^-)$. Therefore we have checked all assumptions of Theorem~\ref{theoCD} with $G = \AH$, since $\AH(0)$ is nilpotent and $p\geq p_1$. Recall from \S \ref{subsec32} that $U(z) z^{\AH(0)}$ is a solution matrix of the differential system   $zZ'(z) = \AH(z)Z(z)$, with $U(z) \in \GLmu(\K[[z]])$, $U(0)=I_\mu$, and $U_{i,j}(z) = \sum_{n=0}^\infty u_{i,j,n} z^n$ with $u_{i,j,n} \in\K$. By~\eqref{eq:CD}, we have 
$$| u_{i,j,n} |_v \leq T(n)^{\mu-1} \mbox{ for any } n\geq 0.$$ 
Here we have used that $p > \mu  $ so that $v_p((\mu-1)!) = v_p (\prod_{j=1}^\mu \binom{\mu}{j}) = 0$.  Since $|D_n|_v = p^{-v_p(D_n)} = p^{-\pe{\frac{\log n}{\log p}}} = T(n)^{-1}$ we obtain:
$$| D_n ^{\mu-1} u_{i,j,n} |_v \leq 1  \mbox{ for any } i,j,n.$$
Now recall that $g(z) = \sum_{n=0}^\infty a_n z^n$, where $a_n$ is given by Eq. \eqref{eqansomme}. Since $p\geq p_3$ we have $|\ell_i P_{k,i}(0)h_{1,j,t}|_v\leq 1$ for any $i$, $j$, $k$, $t$, so that this equation yields 
$$| D_{n+\bz}^{\mu-1} a_n |_v\leq 1 \mbox{ for any } n\geq 0.$$
This upper bound holds for any finite place $v\in\PK$ above any prime number $p\geq p_0$. We shall now deal with primes $p < p_0$.

\subsection{Conclusion of the proof} \label{subsec34}

Since $g$ is a $G$-function, there exists a sequence $(\delta_n)_{n\geq 0}$ of positive  rational  integers such that $\delta_n a_n \in \OK$ and $\delta_n \leq \Delta^{n+1}$ for any $n$,  where $\Delta>1$ is a positive real number that depends only on $g$. Then we let
$$C = \prod_{p< p_0} p^{\pes{\frac{\log \Delta}{\log p}}}$$
where the product ranges over prime numbers $p$ such that $p <p_0$. Now let $p < p_0$, $v\in\PK$ be a place above $p$, and $n\geq 0$. Then we have
$$|C^{n+1}a_n|_v = |a_n|_ v \, p^{-(n+1) \pes{\frac{\log \Delta}{\log p}}} \leq  |a_n|_ v \, p^{-  \frac{\log \delta_n}{\log p}} \leq | \delta_n a_n |_v\leq 1$$
since $\delta_n a_n \in \OK$  and     $|x|_v = p^{-v_p(x)}$ for any $x\in \Z\setminus\{0\}$. Combining this upper bound with the one obtained for $p\geq p_0$, namely $| D_{n+\bz}^{\mu-1} a_n |_v\leq 1 $, we conclude that  $|D_{n+\bz}^{\mu-1} C^{n+1} a_n |_v\leq 1 $ for any finite place $v\in\PK$: for any $n\geq 0$, $
D_{n+\bz}^{\mu-1} C^{n+1} a_n$ is an algebraic integer of $\K$.

\section{Working out an example}\label{sec:example}  \label{sec4}

In this section, we illustrate our results on the example of the generating function of Ap\'ery's numbers 
$a_n=\sum_{k=0}^n \binom{n}{k}^2\binom{n+k}{k}^2$, which appeared in Ap\'ery's famous proof of the irrationality of $\zeta(3)$. The series 
$\mathcal{A}(z)=\sum_{n=0}^\infty a_n z^n$ is a $G$-function because the sequence of integers $(a_n)_{n\ge 0}$ satisfies the linear recurrence
\begin{equation}\label{eq:aperyrec}
(n+2)^3U_{n+2}-(34n^3+153n^2+231n+117)U_{n+1}+(n+1)^3U_n=0,
\end{equation}
which, together with the initial values of $a_n$,  implies that $\mathcal{A}(z)$ is a solution of the differential equation $\mathcal{L}y(z)=0$ with 
$$
\mathcal{L}=z^2(1-34z+z^2)\Big(\frac{d}{dz}\Big)^3 +z(3-153z+6z^2)\Big(\frac{d}{dz}\Big)^2+(1-112z+7z^2) \frac{d}{dz} +z-5.
$$
The singularities of $\mathcal{L}$ are $0$, $( \sqrt{2}-1)^4$ and $(\sqrt{2}+1)^4$. The exponents at $0$ are all equal to $0$, while those at 
$z=(\sqrt{2}\pm 1)^4$ are $0,1,\frac12$.
Moreover, $\mathcal{L}$ is minimal for $\mathcal{A}(z)$, hence it is a $G$-operator. Beukers~\cite{beukers} and Dwork~\cite{dwork} independently observed that 
$\mathcal{L}$ is even a geometric operator. A detailed study was made in~\cite{beukers2};  in particular the following  representation of $\mathcal{A}(z)$ as a period was given, valid for $\vert z\vert <(\sqrt{2}-1)^4$:  
$$
\mathcal{A}(z)=\frac{1}{(2i\pi)^3}\int_{\mathcal{C}} 
\frac{du dv dw}
{1-(1-uv)w-zuvw(1-u)(1-v)(1-w)}
$$
where the algebraic cycle $\mathcal{C}=\{u,v,w \in \mathbb{C}^3 : \vert u\vert=\vert v\vert=\vert w\vert =\frac12\}$; this period is a solution of the Picard-Fuchs equation of the families of $K3$ surfaces $(1-(1-xy)t-zxyt(1-x)(1-y)(1-t)=0)_{z\in \mathbb C}$. See also~\cite{lairez} for further examples of this nature. 
Hence Theorem~\ref{thprinc} applies unconditionally: for any solution $(U_n)_{n\ge 0}$ of \eqref{eq:aperyrec} with $U_n\in \Qbar$ such that $\mathcal{L}(\sum_{n=0}^\infty U_n z^n)=0$, we have that $D_{n+\bz}^{2}C^{n+1} U_n$ are algebraic integers for some positive integers $\bz$ and $C$. Since this applies to the sequence of integers $(a_n)_{n\ge 0}$, we see that the exponent $\mu -1$ of $D_{bn+\bz}$ in Theorem~\ref{thprinc} is not always sharp.

\medskip

For any fixed $\alpha\in \Qbar$, $\alpha\neq 0$, we now consider the local solutions of $\mathcal{L}y(z)=0$ at $z=\alpha$. It amounts to the same to solve locally at $z=0$ the (shifted) differential equation $\mathcal{L}_\alpha y(z)=0$ where 
\begin{multline*}
\mathcal{L}_\alpha=(z+\alpha)^2\big(1-34(z+\alpha)+(z+\alpha)^2\big)\Big(\frac{d}{dz}\Big)^3 
+(z+\alpha)\big(3-153(z+\alpha)+6(z+\alpha)^2\big)\Big(\frac{d}{dz}\Big)^2 \\ 
+\big(1-112(z+\alpha)+7(z+\alpha)^2\big) \frac{d}{dz} +z+\alpha-5.\qquad 
\end{multline*}  

If $\alpha=(\sqrt{2}+1)^4$ or $\alpha=(\sqrt{2}-1)^4$, there is at least one holomorphic solution $\sum_{n=0}^\infty c_n z^n$ of $\mathcal{L}_\alpha y(z)=0$ because at least one exponent at $z=0$ is an integer (they are in fact equal to $0,1,\frac12$). No closed form expression is known for the sequence 
$(c_n)_{n\ge 0}$ but it satisfies the 
linear recurrence (trivially associated to $\mathcal{L}_\alpha$ for any $\alpha\in \mathbb C$)
\begin{multline}\label{eq:aperyrecshift}
(n+1)^3V_n
+
\big((2 \alpha-17) n^2+(6\alpha-51) n+5\alpha-39\big)(2n+3)
V_{n+1}
\\+
\big((6\alpha^2 -102\alpha+1)n^2+(24\alpha^2-408\alpha+4)n+25\alpha^2-418\alpha+4\big)(n+2)
V_{n+2}
\\+\alpha(2\alpha^2-51\alpha+1)(2n+5)(n+3)(n+2)V_{n+3}
\\+
\alpha^2(\alpha^2-34\alpha+1)(n+4)(n+3)(n+2)
V_{n+4}=0.
\end{multline} 
Note that for $\alpha=(\sqrt{2}\pm 1)^4$, the recurrence is of order $3$ only because $\alpha^2-34\alpha+1=0$. 
Since $\mathcal{L}_\alpha$ is a geometric operator, we can apply Theorem~\ref{thprinc} unconditionally to any solution $(V_n)_{n\ge 0}$ of \eqref{eq:aperyrecshift} with $V_n\in \Qbar$ such that $\mathcal{L}_\alpha(\sum_{n=0}^\infty V_n z^n)=0$, which includes $(c_n)_{n\ge 0}$. It follows that there exist some positive integers $\bz$ and $C$ such that 
$D_{2n+\bz}^2C^{n+1} V_n$ is an algebraic integer for all $n\ge 0$.

We now deal with the case $\alpha \notin\{0, (\sqrt{2}+1)^4, (\sqrt{2}-1)^4\}$. Then $z=0$ is a regular point of $\mathcal{L}_{\alpha}$ and thus its exponents are integers. In fact they are equal to $0,1$ and $2$ for any such $\alpha$. Let us denote by $\sum_{n=0}^\infty c_{j,n}z^n\in \mathbb Q[[z]]$, $j=1,2,3$, a basis of local holomorphic solutions at $z=0$ of $\mathcal{L}_\alpha y(z)=0$. They are $G$-functions but closed forms for the coefficients $c_{j,n}$ are not known.
Since $\sum_{n=0}^\infty c_{j,n}z^n$ are independent solutions of $\mathcal{L}_\alpha y(z)=0$, the three sequences 
$(c_{j,n})_{n\ge 0}$ are independent solutions of the recurrence~\eqref{eq:aperyrecshift}. 
The latter is of order $4$ so that the $(c_{j,n})_{n\ge 0}$ do not generate all the solutions $(V_n)_{n\ge 0}$ of \eqref{eq:aperyrecshift}. However, they generate all the sequences $(V_n)_{n\ge 0}$ with $V_n\in \Qbar$ such that $\mathcal{L}_\alpha(\sum_{n=0}^\infty V_n z^n)=0$. By Theorem~\ref{thprinc}, we are ensured unconditionally that for any such $(V_n)_{n\ge 0}$, we have $D_{n+\bz}^{2}C^{n+1} V_{n}\in \mathbb Z$ for all $n\ge0$ and for some positive integers $\bz$,~$C$.

\medskip

If $(U_n)_{n\ge 0}$ is a solution of the recurrence  \eqref{eq:aperyrec} with $U_n\in \Qbar$ but $\mathcal{L}(\sum_{n=0}^\infty U_n z^n)\neq 0$, then in fact $\mathcal{L}(\sum_{n=0}^\infty U_n z^n)=P(z)$ for some polynomial $P(z)\in \Qbar[z]$ of degree $k$ say. Hence, $\sum_{n=0}^\infty U_n z^n$ is a $G$-function because $(\frac{d}{dz})^{k+1}\circ \mathcal{L}$ is still a $G$-operator. It is of order $k+4$, of geometric type, with integer exponents at $z=0$. Hence, we can apply Theorem~\ref{thprinc} unconditionally: there exist some positive integers $\bz, C$ such that $D_{n+\bz}^{k+3}C^{n+1} U_n$ are algebraic integers for all $n\ge0$. This situation applies to the second sequence of rational numbers used by Ap\'ery to prove the irrationality of $\zeta(3)$, namely
$$
\widehat{a}_n=\sum_{k=0}^n \binom{n}{k}^2\binom{n+k}{k}^2\bigg(\sum_{m=1}^n \frac{1}{m^3}+\sum_{m=1}^k \frac{(-1)^{m-1}}{2m^3\binom{n}{m} \binom{n+m}{n}}\bigg).
$$
Indeed $(\widehat{a}_n)_{n\ge 0}$ is a solution of the recurrence  \eqref{eq:aperyrec} but 
$\mathcal{L}(\sum_{n=0}^\infty \widehat{a}_n z^n)=5$ (see~\cite{beukers2}). We thus have $D_{n+\bz}^3 C^{n+1}\widehat{a}_n\in \mathbb Z$ for some integer $\bz\ge 0$. In this case, this is sharp as it is known that the exponent $3$ can not be improved; moreover we can take $\bz=0$ and $C=1$. Similar observations can be formulated for the solutions $(V_n)_{n\ge 0}$ of the recurrence \eqref{eq:aperyrecshift} such that $V_n\in \Qbar$ but $\mathcal{L}_{\alpha}(\sum_{n=0}^\infty V_nz^n)\neq 0$.

\providecommand{\bysame}{\leavevmode ---\ }
\providecommand{\og}{``}
\providecommand{\fg}{''}
\providecommand{\smfandname}{\&}
\providecommand{\smfedsname}{eds.}
\providecommand{\smfedname}{ed.}
\providecommand{\smfmastersthesisname}{M\'emoire}
\providecommand{\smfphdthesisname}{Th\`ese}

\bigskip

\noindent S. Fischler, 
Laboratoire de Math\'ematiques d'Orsay, Univ. Paris-Sud, CNRS, Universit\'e Paris-Saclay, 91405 Orsay, France.

\medskip

\noindent T. Rivoal, Institut Fourier, CNRS et Universit\'e Grenoble Alpes, CS 40700, 
 38058 Grenoble cedex 9, France

\bigskip

\noindent Keywords: $G$-functions, $G$-operators, $p$-adic differential equations, Christol-Dwork Theorem.

\medskip

\noindent MSC 2000: 11J72, 12H25.

\end{document}